\documentclass{article}%
\usepackage{amsmath}%
\usepackage{amsfonts}%
\usepackage{amssymb}%
\usepackage{graphicx}%

\def\Bbb{\mathbb}
\newcommand{\sign}{\text{sign}}

\def\l{\lambda }
\def\a{\alpha }
\def\b{\beta }

\def\s{\sigma}

\def\per{\mathop{\rm per}}

\def\T{\mathop{\footnotesize{t}}}

\begin{document}
\title{An update on a few  permanent conjectures
\thanks{In memory of Professors Marvin Marcus (1927-2016) and Ingram Olkin (1924-2016). Presented at the 5th International Conference on Matrix Analysis and Applications, Fort Lauderdale, Florida, USA, Dec. 17-20, 2015.}
}
\author{Fuzhen Zhang\thanks{E-mail: zhang@nova.edu.}\\
\footnotesize{Department of Mathematics, Halmos College,
Nova Southeastern University}\\
\footnotesize{3301 College Avenue, Fort Lauderdale, FL 33314, USA}\\
}

\date{}
\maketitle

{\small \bf ABSTRACT:} We review and update on a few conjectures concerning matrix permanent that are easily stated, understood, and accessible to general math audience. They are:
 Soules permanent-on-top conjecture${}^\dagger$, Lieb permanent dominance conjecture, Bapat and Sunder conjecture${}^\dagger$ on Hadamard product and diagonal entries, Chollet conjecture
  on Hadamard product, Marcus conjecture on permanent of permanents, and several other conjectures. Some of these conjectures are recently settled; some are still open.
We also raise a few new questions for future study. (${}^\dagger$conjectures have been recently settled negatively.)


\section{Introduction}
Computed from the elements of a square matrix, the {determinant} of a square matrix is one of the most useful and important concepts in mathematics.
Of many matrix functions, {\em permanent} is another important one. It arises naturally in the study of the symmetric tensors in multilinear algebra (see, e.g., \cite{MerMul97}); it also plays a role in combinatorics
(see, e.g., \cite{Lint01}).
Both terms were introduced in the 1800s (see, e.g., \cite[p.\,1]{Minc78}); and they are still useful  in research.

One of the intriguing problems on permanent is the so-called van der Waerden conjecture (theorem) regarding the minimum value of the permanent on the Birkhoff polytope of doubly stochastic matrices. It was conjectured by van der Waerden in 1926 and resolved
in the affirmative by Egory\^{c}ev and Falikman independently in 1981 (see, e.g., \cite{Minc83}). Second to the van der Waerden conjecture is, in my opinion, the {\em permanence dominance conjecture}
(see below). Since its appearance in the mid-1960s,
it has drawn much attention of the mathematicians in the area; it remains open as one of the most important unsolved problems
in linear algebra and matrix theory. A closely related and stronger statement
is the {\em permanent-on-top} conjecture which appeared about the same time and is recently shown to be false with the help of computation utility.

We begin by the permanent-on-top (POT) conjecture and  two  chains of conjectures on permanent starting with the POT conjecture,
reviewing and updating their developments and status. We also include a few other conjectures
on permanent of our interest that are easily stated but remain unsolved or even have no progress over the decades.
A few comprehensive surveys on permanent, including many open conjectures and problems,  are the Minc's monograph {\em Permanents} 1978 \cite{Minc78}, and subsequent articles {\em Theory of Permanents 1978--1981} \cite{Minc83}
and {\em Theory of Permanents 1982--1985} \cite{Minc87}, followed by
Cheon and Wanless' {\em An update on Minc's survey of open problems involving permanents} 2005 \cite{CW05}.

\section{Preliminaries}
Let $M_n$ be the set of all $n\times n$ complex matrices. We assume $n\geq 2$ in the paper.  For an $n\times n$ matrix $A=(a_{ij})$, the determinant
$\det A$ of $A$ is $\sum_{\sigma \in S_{n}}
\prod _{t=1}^{n}  \sign (\sigma) a_{t\sigma (t)},$ where $S_{n}$ is
the symmetric group of degree $n$,
while the permanent of $A=(a_{ij})$,
denoted by $\per A$, is  simply defined as
$\per A=\sum_{\sigma \in S_{n}}
\prod _{t=1}^{n}a_{t\sigma (t)}$.
For $A\in M_n$, by writing $A>0$ (resp. $A\geq 0$) we mean that $A$ is  positive (resp. semi) definite,
that is, $x^*Ax>0$ (resp. $x^*Ax\geq 0$) for
all  column vectors $x\in \Bbb C^n$, where $x^*=(\bar{x})^t$ is the
conjugate transpose of $x$.
It is known that if $A$ is positive semidefinite, then $0\leq \det A \leq \prod_{i=1}^na_{ii} \leq \per A$. Determinant and permanent may be regarded as sister functions of positive semidefinite matrices as many of their results   exist side by side. For example, when $A$ and $B$ are $n\times n$ and positive semidefinite, then
$\det (A+B)\geq \det A+\det B$ and $\per (A+B)\geq \per A+\per B$.
(Note:  the directions of inequalities for ``det" and ``per" sometimes are the same and some other times are reversed.) A permanent analogue of a determinant result does  not always exist. For instance,
if $A\geq 0$ is written as $A=R+Si$, where $R$ and $S$ are Hermitian, then
$\det A\leq \det R$ (the Robertson-Taussky inequality); but $\per A$ and $\per R$ are incomparable in general \cite{JZ95}.

Determinant and permanent are special {\em generalized matrix functions}.
Let $S_n$ be the symmetric group of degree $n$,
 $H$ be a subgroup of $S_n$,  and $\chi$ be a character
 on $H$. The generalized matrix function of an $n\times n$ matrix
$A=(a_{ij})$ with respect to $H$ and $\chi$ is
defined by (see, e.g., \cite[p.~124]{MarcusMul1})
$$d^H_{\chi}(A)=\sum _{\s \in H}\chi (\s)\prod _{t=1}^n a_{t\s (t)}$$

Setting $H=S_n$ and $\chi (g)=\sign (g)=\pm 1$ according to  $g\in H$ being even or odd, we have the determinant $\det A$;
putting $H=S_n$ and $\chi (g)=1$ for all $g$,  we get the permanent $\per A=\sum_{\s \in S_n}\prod _{i=1}^n a_{i\s (i)}$.
The product of the main diagonal entries of $A$, $h(A)=a_{11}a_{22}\cdots a_{nn}$,   known as the
Hadamard matrix function, is also a generalized matrix function by taking $H=\{e\}$, where $e$ is the group identity of $S_n$.
Specially, if  $H=S_n$ and $\chi$ is an irreducible character of $S_n$, then the permanent
$d^H_{\chi}(A)$ is  referred to as  {\em immanant}.

Let $A=(a_{ij})$ be a square matrix. We write $A=(A_{ij})\in M_n(M_m)$ if $A$ is an
$n\times n$ partitioned matrix in which each block $A_{ij}$ is $m\times m$.
(Thus,  $A$ is an $mn\times mn$ matrix.) In case $n=2$, we  write explicitly
$A=\left ( {A_{11} \atop  A_{21}} {A_{12}\atop A_{22}} \right ).$

 Listed below are a few well-known milestone inequalities for determinant, permanent, and generalized matrix functions of positive semidefinite matrices. They stand out because of their elegancy in form and importance in applications.  Let $A=(a_{ij})$ be an $n\times n$ positive semidefinite matrix with a partitioned form $A=\left ( {A_{11} \atop  A_{21}} {A_{12}\atop A_{22}} \right )$, where $A_{11}$ and $A_{22}$ are square submatrices (possibly of different sizes) of $A$. We have the following inequalities concerning $A$.
\bigskip

Hadamard inequality (1893):
\begin{equation}\label{Had}
\det A\leq \prod_{i=1}^n a_{ii}
\end{equation}

Fischer inequality (1908):
\begin{equation}\label{Fis}
\det A\leq \det A_{11} \det A_{22}
\end{equation}

Schur inequality (1918):
\begin{equation}\label{Sch}
\det A\leq \frac{1}{\chi (e)} d_{\chi}^H(A)
\end{equation}

Marcus inequality (1963):
\begin{equation}\label{Mar}
\per A \geq  \prod_{i=1}^n a_{ii}
\end{equation}

Lieb inequality (1966): %
\begin{equation}\label{Lie}
  \per A\geq \per A_{11} \per A_{22}
\end{equation}

Schur inequality (\ref{Sch}) implies Fischer inequality (\ref{Fis}) which implies Hadamard inequality (\ref{Had}). 
(\ref{Mar}) and (\ref{Lie}) are permanent analogues of (\ref{Had}) and (\ref{Fis}), respectively.
Inequality (\ref{Lie}) (in a more general form)  was first conjectured by Marcus and Newman in 1965 \cite{MM65Per} and proved by Lieb in 1966 \cite{Lieb66}, and re-proved by
Djokovi\v{c} in 1969 \cite{Djo69}.

If we compare (\ref{Had}) and (\ref{Mar}), and (\ref{Fis}) and (\ref{Lie}),
in view of (\ref{Sch}), we can naturally ask if the right hand of
(\ref{Sch}) is bounded by $\per A$, that is, if the following holds:
\begin{equation}\label{LieConj}
 \frac{1}{\chi (e)} d_{\chi}^H(A)\leq \per A
\end{equation}

This is a conjecture of Lieb  also known as the {\em permanent dominance conjecture}. 
It was stated explicitly by Lieb in 1966 \cite{Lieb66}
when studying a similar problem of Marcus and Minc 1965 \cite{MM65Per}:
 Under what conditions on $\chi $ and $H$ will the following inequality
 hold for all $n\times n$ positive semidefinite $A$
\begin{equation}\label{Problem2a}
\sum_{\s \in H}\chi (\s) \prod_{i=1}^na_{i\s (i)}\leq \per A?
\end{equation}

(\ref{LieConj}) and (\ref{Problem2a}) are the same when  $\chi$ is degree 1, i.e., $\chi(e)=1$.

A great amount of  effort was made in attempting to solve the Lieb permanent dominance  conjecture, especially from the late 1960s to early 1990s.
Motivated by Lieb conjecture \cite{Lieb66}, Soules proposed  in his 1966 Ph.D. thesis
(see also Merris 1987 \cite{Mer87}) 
 a  conjecture stronger than the permanent dominance conjecture.
  Soules conjecture
 states that the permanent of a positive semidefinite matrix $A$ is the largest eigenvalue of the Schur power matrix of $A$ (see Section \ref{Sec3}).
 This is referred to as {permanent-on-top} (POT) conjecture.
 In fact, POT conjecture is the strongest among several permanent conjectures (see, e.g., \cite{CW05}):
$$\begin{array}{lll}
      &   &     \mbox{Lieb per-dom}  \longrightarrow \mbox{Marcus per-in-per}       \\
      & \nearrow  &     \\
\mbox{Soules POT${}^\dagger$}     &   & \\
      & \searrow  &  \\
       &   & \mbox{Bapat \& Sunder${}^\dagger$}  \longrightarrow \mbox{Chollet}
\end{array}$$

It is known now that Soules  and Bapat \& Sunder conjectures are false; but the other three, the Lieb, Marcus and Chollet conjectures, are still open.
In view of recent developments in the area and with more advanced computation tools available nowadays, it is worth revisiting these long standing conjectures.



\section{The permanent-on-top conjecture is false}\label{Sec3}

In a recent publication, Shchesnovich 2016 \cite{Sch16} settled in the negative a long standing conjecture on permanent; it shows that the
 permanent-on-top (POT) conjecture is false. The POT conjecture was originally formulated by Soules in his Ph.\,D. dissertation 1966 \cite[Conjecture, p.\,3]{Sou66} (formally published in Minc's 1983
 \cite[p.\,249]{Minc83}) in attempting to answer the Lieb conjecture \cite{Lieb66}.
 It
 stated that the permanent of a positive semidefinite matrix $A$ was the largest eigenvalue of the {Schur power matrix} of $A$.
 For an $n\times n$ matrix $A=(a_{ij})$, the {\em Schur power matrix} of $A$, denoted by
 $\pi(A)$, is an $n!\times n!$ matrix with entries $s_{\a \b}$ lexicographically indexed by permutations $\a, \b\in S_n$:
 $$s_{\a \b}=\prod_{t=1}^n a_{\a (t) \b (t)}$$

 As a principal submatrix of the $n$-fold tensor (Kronecker) power $\otimes^n A$, the Schur power matrix ${\pi (A)}$ inherits many properties of $A$. For instance,
 if $A$ is positive semidefinite, then so is its Schur power; and the eigenvalues of ${\pi (A)}$ interlace those of $\otimes^n A$. Observing that for each fixed $\alpha \in S_n$,
 $$\sum_{\beta\in S_n} s_{\a \b}=\sum_{\beta\in S_n} \prod_{t=1}^n a_{\a (t) \b (t)}= \per A$$

 We see that every row sum of ${\pi (A)}$ is $\per A$ \cite{Sou94}. Consequently,
 $\per A$ is an eigenvalue of ${\pi (A)}$. It is also known (and not difficult to show) that $\det A$ is an eigenvalue of ${\pi (A)}$ too. In fact, $\det A$ is the smallest eigenvalue of ${\pi (A)}$ when $A\geq 0$; see \cite[p.\,221]{MerMul97}.
  For positive semidefinite $A$,  the POT conjecture claimed that $\per A$ was the largest eigenvalue of $\pi(A)$. Note that POT conjecture is independent of subgroup $H$ and character $\chi$ in (\ref{LieConj}).

 The counterpart of  the POT  conjecture  for determinant is that  $\det A$ is the smallest eigenvalue of $S(A)$ which is shown implicitly in  Schur 1918 \cite{ISchur18},
 as  pointed out
 by Bapat and Sunder 1986 \cite[p.\,154]{BS86}.
 Bapat and Sunder 1986 \cite{BS86} proved that the  POT conjecture is true for $n\leq 3$ (and also provided an equivalent form of the POT conjecture). 
 POT conjecture involved no subgroup and character and
  it would yield Lieb conjecture (in the case of degree 1 characters as assumed in \cite[p.\,3]{Sou66} and more generally in
 Merris 1987 \cite{Mer87}) 
 because the value on the left hand side of
 (\ref{LieConj}) is contained in the numerical range of the Schur power matrix.
 If $W(X)$ denotes the numerical range of a square matrix $X$, then the POT  conjecture is equivalent to the statement that for $A\geq 0$,
 \begin{equation}\label{rep}
 W\left  (\pi(A)\right ) =[\det A,\; \per A]
 \end{equation}

 In fact, Soules POT conjecture was the strongest among several permanent conjectures; see the conjecture chains in the previous section; see also, e.g., \cite{CW05}.

 Notice that
 if $A$ is an $n\times n$ positive semidefinite matrix, then $x^*Ax\leq \l_{\max}(A) x^*x$ for all column vectors $x\in \Bbb C^n$,
 where $\l_{\max}(A)$ is the largest eigenvalue of $A$.
 Shchesnovich 2016 \cite{Sch16} presented an example
  of $5\times 5$ positive semidefinite matrix $H$ (thus $\pi(H)$ is $5!\times 5!$, i.e., $120\times 120$) and  a column vector $X$ of $5!=120$  components
   (which is omitted here as it is too large) such that
 $$X^*\pi(H)X> \per (H) X^*X$$
 where
 $$H= \left ( \begin{array}{ccccc}
 40  & 14 -22i  &  -4 - 8i &  8 +16i & 22 -36i\\
  14 +22i & 22  & 16 -14i   & -9 - i & 23 -12i\\
  -4 + 8i &16 +14i & 52     &  8 -34i & 14 -30i\\
   8-16i & -9 + i &  8 +34i   & 34  & 18 -19i\\
  22 +36i & 23 +12i & 14 +30i   & 18 +19i & 75
 \end{array} \right )=u^*u+v^*v$$ with
$$u=(4+2i, 2-3i, -4-4i, -3+4i, 1),\;   v=(2-4i, -3i, 2-4i, -3i, -5-7i).$$
$H$ is a $5\times 5$  positive semidefinite matrix of rank 2 having eigenvalues 0, 0, 0, 91, and 132, while its Schur power matrix $\pi(H)$ is of rank 27.

 The Schur power matrix $S(H)$ is of order $5!=120$. It would be difficult to compute and display $\pi(H)$. A computation utility
 is needed.
 It is also demonstrated in \cite{Sch16} through characteristic polynomial  that $\l_{\max}(\pi(H))>\per H$:
$$\l_{\max} ({\pi}(H))=320\big (2185775+\sqrt{160600333345}\,\big )  > \per (H)=814016640$$

For an $n$-square matrix $A$, the Schur power matrix $\pi(A)$ as a principal submatrix of the tensor power $\otimes ^n A$ may have interest in its own right.
It is easy  to show that $\{w_1w_2\cdots w_n \mid w_i\in W(A), i=1, 2, \dots, n\}\subseteq W(\otimes^n A)$.

\begin{enumerate}

\item[]{\bf Question 1:} Can one find a counterexample of size $n=4$ for the POT conjecture?
 As is known, the POT conjecture is true for $n\leq 3$. With a counterexample
 of $5\times 5$, Shchesnovich 2016 \cite[p.\,198]{Sch16} states that a $4\times 4$ counterexample was not found despite an extensive search.

\item[]{\bf Question 2:} Referring (\ref{rep}), what is exactly the maximum value (which may exceed $\per A$)  in the numerical range of the Schur power matrix of a positive semidefinite matrix? Of a Hermitian matrix? That is, find
    $\l_{\max} (\pi(A))$. This question is especially of interest as
    the numerical range of a Hermitian matrix is a closed interval $[d, p]$.
    What are the left and right end points, i.e, $d$ and $p$, of the interval?
    A result of Pate 1989 \cite[Theorem 4]{Pate89} 
    implies  that $\langle \pi(A) x, x\rangle\leq \per A$ for all unit vectors $x$ in {certain} (but not all) subspaces of $\Bbb C^{n!}$.

\item[]{\bf Question 3:} Investigate the numerical range of the Schur power matrix $\pi(A)$ for a general $n$-square matrix $A$. It is  contained in
    the numerical range of $\otimes^n A$ because $\pi(A)$ is a principal submatrix of
    the latter.

 \end{enumerate}

Note: There has been a confusion about the names of the conjectures in the literature. In 1987 \cite[Conjecture 1]{Mer87}, Merris called inequality (\ref{LieConj})
the {\em Permanental Dominance Conjecture}; in the same book, Johnson 1987 \cite[p.\,168]{Joh87}
 called (\ref{LieConj}) the {\em permanent-on-top conjecture}, while in their introductions,
Bapat \cite{Bapat07} and Zhang \cite{{FZ13Per}} considered the two 
being the same.
Besides, there have been several variations of the  POT conjecture;
see Merris 1987 \cite{Mer87}.

Now it is known that the Soules POT conjecture is false. The weaker one,  Lieb permanence dominance conjecture, is still open. It has been proven true for special cases of $H=S_n$ and some characters $\chi$ (see the next section).

\section{The permanent dominance conjecture is open}

In 1965, Marcus and Minc proposed the following question regarding the permanent of  positive semidefinite matrices with a subgroup of the permutation group and group character 
(see also 
\cite[p.\,158]{Minc78},  \cite[p.\,244]{Minc83}, and \cite[p.\,132]{Minc87}):

{Problem 2 of \cite[ p.\,588]{MM65Per}.  \em
Let $H$ be a subgroup of $S_n$ and let $\chi$ be a character of degree 1 of $H$. Under what conditions on $\chi $ and $H$ will the inequality
\begin{equation}\label{Problem2}
\sum_{\s \in H}\chi (\s) \prod_{i=1}^na_{i\s (i)}\leq \per A
\end{equation}
 hold for all $n\times n$ positive semidefinite Hermitian $A$?}

In 1966, Lieb stated that the answer to the above question is given by the conjecture below (referred to as Conjecture $\alpha$ in Lieb 1966 \cite{Lieb66}):
\medskip

\noindent
{\bf Lieb permanent dominance conjecture 1966 \cite{Lieb66}.}
{\em Let $H$ be a subgroup of the symmetric group $S_n$ and let $\chi$ be a character of degree $m$ of $H$. Then
\begin{equation}\label{Lieb2}
\frac{1}{m}\sum_{\s \in H}\chi (\s) \prod_{i=1}^na_{i\s (i)}\leq \per A
\end{equation}
 holds for all $n\times n$ positive semidefinite Hermitian $A$}.

The permanent dominance conjecture (\ref{Lieb2}) claims that $\per A$ is the smallest (best) upper-bound of $\frac{1}{m}\sum_{\s \in H}\chi (\s) \prod_{i=1}^na_{i\s (i)}$ for all subgroups and characters.
Merris 1983 \cite{Mer83Isr} gave an upper bound
\begin{equation}\label{Merupper}
\frac{1}{m}\sum_{\s \in H}\chi (\s) \prod_{i=1}^na_{i\s (i)} \leq
\left ( \prod_{i=1}^n (A^n)_{ii}\right )^{1/n}\leq \frac{1}{n} \sum_{i=1}^n \l_i^n
\end{equation}
where $(A^n)_{ii}$ denotes the $i$th diagonal entry of $A^n$ and the $\l_i$s are the eigenvalues of $A$. The middle term in (\ref{Merupper})
 is no less than $\per A$ as  noted in \cite[p.\,215]{Mer87}.


 Historically, in his Ph.D. dissertation, Soules 1966 \cite[p.\,3]{Sou94} (where characters were assumed to have degree 1) proposed the later so-called POT conjecture, stating that
 it had  the permanent dominance inequality as its consequence. This method for approaching
  (\ref{Problem2}) and (\ref{Lieb2}) was announced by Soules
at the 1981 Santa Barbara Multilinear Algebra Conference (see \cite[p.\,245]{Minc83}).

In 1986 at the Auburn Matrix Conference, Merris
 called Lieb's conjecture {\em Permanent Dominance Conjecture}
 (see \cite[p.\,135]{Mer86}, and also \cite{Mer87}). This term has later been used by Minc \cite[p.\,135]{Minc87}, Soules \cite{Sou94}, Pate \cite{Pate99},
 Cheon and Wanless \cite{CW05}, and many others.
 So it has become known as the {\em permanent dominance conjecture}.
Merris \cite[p.\,216]{Mer87} 
 explicitly
 pointed out that the Soules POT conjecture implied the Lieb permanent dominance  conjecture. In his 1994 paper \cite{Sou94}, Soules gave a clear and brief explanation of this.
 In fact, Souels POT conjecture would imply several permanent conjectures (see, e.g., \cite{CW05}).
A great amount of work was devoted to the permanent dominance conjecture since its appearance to early 1990s. As an open problem on permanent, this conjecture is arguably most interesting and intriguing; it involves group theory, combinatorics, and matrix theory, of course.

Several paths have been taken to attack or approach the Lieb conjecture. Observe that
the right hand side of (\ref{Lieb2})  is independent of the choices of the subgroup and the character. Immediate cases to be studied are:
$H=S_n$, $\chi$ is the principal character (i.e., $\chi (e)=1)$, or $\chi$ is irreducible.

Great amount work has been done for immanant
(with $H=S_n$ and irreducible $\chi$) and some significant results
have been obtained. It is known \cite{Pate98, Pate99} that
the permanent dominance conjecture is true
for immanants with $n\leq 13$.

The interest and attention  on this interesting and rather difficult conjecture appeared fading away in the last decades.
In the author's opinion, there is still a long way to go for a complete solution of the general conjecture.  Cheon and Wanless  2005 \cite[p.\,317]{CW05} set an excellent account in detail about the permanent dominance conjecture. Other references on this conjecture are
Merris 1987 \cite{Mer87}, Johnson 1987 \cite{Joh87}, and references therein.


\section{Bapat \& Sunder conjecture is false}

Shchesnovich 2016 \cite{Sch16}
devotes a section discussing a permanent conjecture of Bapat and Sunder 1985 \cite{BS85} that is weaker than the POT conjecture and explains that the weaker one may be true for a physical
reason (quantum theory). The counterexample provided in \cite{Sch16}  settles  POT conjecture negatively but does not work for the Bapat \& Sunder conjecture.

Recall the
Oppenheim inequality \cite{Opp30} for positive semidefinite matrices: Let $A=(a_{ij})$ and $B=(b_{ij})$ be $n\times n$ positive semidefinite matrices. Then
\begin{equation}\label{Opp}
\det (A\circ B) \geq \det A \prod_{i=1}^nb_{ii}\geq \det A \det B
\end{equation}
where $A\circ B=(a_{ij}b_{ij})$ is the Hadamard (a.k.a Schur) product of $A$ and $B$.
In view of (\ref{Opp}) for permanent, Bapat and Sunder 1985 \cite[p.\,117]{BS85}
and 1986 \cite[Conjecture 2]{BS86} raised the
question (conjecture) whether 
\begin{equation}\label{BS85}
\per (A\circ B)\leq \per A\; b_{11}\cdots b_{nn}
\end{equation}

It is not difficult to reduce Bapat \& Sunder conjecture to the case of correlation matrices (see, e.g.,  \cite{ZFZ89}).
By a {\em correlation matrix},  we mean  a positive semidefinite matrix all whose main
diagonal entries are equal to 1. That is, Bapat \& Sunder conjecture is equivalent to the statement that
{
if $A$ and $B$ are correlation matrices, then
\begin{equation}\label{Z89}
\per(A\circ B)\leq \per A
\end{equation}
}

Bapat \& Sunder conjecture is weaker than the permanent-on-top conjecture
(which is now known false); see Merris 1987 \cite[Conjecture 10, p.\,220]{Mer87}.
 Drury 2016 \cite{Dru16} found a  counterexample and
settled this negatively. Let $A$ be
{\small $$
\begin{pmatrix}
1&0&\frac{1}{\surd 2}&\frac{1}{\surd 2}&\frac{1}{\surd 2}&\frac{1}{\surd 2}&\frac{1}{\surd 2}\\
0&1&\frac{1}{\surd 2} e^{\frac{4}{5} i \pi}&\frac{1}{\surd 2} e^{\frac{2}{5} i \pi}&\frac{1}{\surd 2}&\frac{1}{\surd 2} e^{-\frac{2}{5} i \pi}&\frac{1}{\surd 2} e^{-\frac{4}{5} i \pi}\\
\frac{1}{\surd 2}&\frac{1}{\surd 2} e^{-\frac{4}{5} i \pi}&1&\cos\left(\frac{1}{5} \pi\right) e^{-\frac{1}{5} i \pi}&\cos\left(\frac{2}{5} \pi\right) e^{-\frac{2}{5} i \pi}&\cos\left(\frac{2}{5} \pi\right) e^{\frac{2}{5} i \pi}&\cos\left(\frac{1}{5} \pi\right) e^{\frac{1}{5} i \pi}\\
\frac{1}{\surd 2}&\frac{1}{\surd 2} e^{-\frac{2}{5} i \pi}&\cos\left(\frac{1}{5} \pi\right) e^{\frac{1}{5} i \pi}&1&\cos\left(\frac{1}{5} \pi\right) e^{-\frac{1}{5} i \pi}&\cos\left(\frac{2}{5} \pi\right) e^{-\frac{2}{5} i \pi}&\cos\left(\frac{2}{5} \pi\right) e^{\frac{2}{5} i \pi}\\
\frac{1}{\surd 2}&\frac{1}{\surd 2}&\cos\left(\frac{2}{5} \pi\right) e^{\frac{2}{5} i \pi}&\cos\left(\frac{1}{5} \pi\right) e^{\frac{1}{5} i \pi}&1&\cos\left(\frac{1}{5} \pi\right) e^{-\frac{1}{5} i \pi}&\cos\left(\frac{2}{5} \pi\right) e^{-\frac{2}{5} i \pi}\\
\frac{1}{\surd 2}&\frac{1}{\surd 2} e^{\frac{2}{5} i \pi}&\cos\left(\frac{2}{5} \pi\right) e^{-\frac{2}{5} i \pi}&\cos\left(\frac{2}{5} \pi\right) e^{\frac{2}{5} i \pi}&\cos\left(\frac{1}{5} \pi\right) e^{\frac{1}{5} i \pi}&1&\cos\left(\frac{1}{5} \pi\right) e^{-\frac{1}{5} i \pi}\\
\frac{1}{\surd 2}&\frac{1}{\surd 2} e^{\frac{4}{5} i \pi}&\cos\left(\frac{1}{5} \pi\right) e^{-\frac{1}{5} i \pi}&\cos\left(\frac{2}{5} \pi\right) e^{-\frac{2}{5} i \pi}&\cos\left(\frac{2}{5} \pi\right) e^{\frac{2}{5} i \pi}&\cos\left(\frac{1}{5} \pi\right) e^{\frac{1}{5} i \pi}&1
\end{pmatrix}
$$}

Drury's matrices $A$ (of rank 2 with complex entries)  and $B=A^{\T}=\bar{A}$ are $7\times 7$ correlation matrices
with $\per A=45$ and $\per (A\circ B)=\frac{6185}{128}$,
 having the ratio
$\per (A\circ B)/\big (\per A \prod_{i=1}^7 b_{ii}\big )=1237/1152>1$.
Moreover, $\l_{\max}(A)=3.5$, while $\l_{\max}(\pi(A))=525/8 = 65.625$.
Very recently, Drury \cite{Dru160519}
    provided a {\em real} counterexample of size $16\times 16$ with rank 3. So $n=16$ is the smallest
    size of known real positive semidefinite matrices that fail the POT conjecture.
    Note that Drury's counterexamples also disprove the POT conjecture as the latter is stronger.

\begin{itemize}
\item[]{\bf Question 4}: Can one find a counterexample of size $n<7$, say $n=4$, with complex (or real) entries,
 for the Bapat \& Sunder conjecture?
 \medskip

\item[]{\bf Question 5:} Can one find a {\em real} counterexample of size $n< 16?$
    That is to say, does the Bapat \& Sunder conjecture (and also the POT conjecture)  hold true for {\em real} positive semidefinite matrices of sizes no more than 15?
    Shchesnovich 2016 \cite[p.\,198]{Sch16} states that all found counterexamples are complex and rank deficient.  Drury's  counterexamples are also rank deficient. Soules 1994 \cite{Sou94} showed certain necessary conditions for a {\em real} positive semidefinite matrix to fail the conjecture. (Note: the statement ``By Theorem 1, the conjecture can only fail at a singular matrix" on page 222, line 19 of \cite{Sou94} is incorrect.) In particular, for the case of real and $n=4$, if POT fails, then there exists a {\em singular} matrix for which the POT is false. It follows that the POT conjecture would be true for $n=4$ if one could prove that it is true for all {\em singular} positive semidefinite matrices of size $n=4$.

\end{itemize}

Denote by ${\cal C}_{n}$   the collection of all
$n\times n$
  complex   correlation matrices.
The set ${\cal C}_{n}$   can be thought of as a
subset of $\Bbb C^{n^{2}}$ and it is compact and convex. The compactness of ${\cal C}_{n}$ follows from
the fact that $A\geq 0$ if and only if
all principal submatrices of $A$ have nonnegative determinants.
We also see  that ${\cal C}_{n}$ is closed
under the Hadamard product.
Given $A\in {\cal C}_{n}$, we define a function on ${\cal C}_{n}$ by
\begin{equation}\label{f(X)}
f_A(X) =\per(A\circ X), \;\;  X \in {\cal C}_{n}
\end{equation}

Since ${\cal C}_{n}$ is compact and $f_A$ is continuous, there exists
a correlation matrix depending on $A$, referred to as {\em maximizer} of $A$ and denoted by $M_{A}$,  such that
$$f(M_{A})=\max_{X\in {\cal
C}_{n}}  \per(A\circ X)$$ i.e., the maximal value can be attained.
Drury's example shows that it is possible that
a (irreducible) maximizer has all entries with moduli  less than 1.
Several properties of maximizing matrices are presented in Zhang 2013 \cite{FZ13Per}.

In view of (\ref{Opp}), we see that the determinant is a ``Hadamard-dilation" function on
${\cal C}_{n}$ in the sense that $\det A\leq \det (A\circ X)$ for $A\in {\cal C}_{n}$ and
any $X\in {\cal C}_{n}$. Equivalently,
$\max\{\det A, \det B\} \leq \det (A\circ B)$ for $A, B\in {\cal C}_{n}$.
 In contrast, the permanent is no ``Hadamard-compression" on ${\cal C}_{n}$.
\medskip

\begin{figure}[h] 
  \centering
  \includegraphics[width=1.3in]{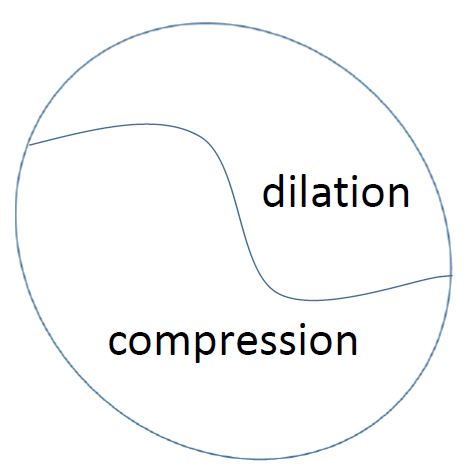}
  \caption{det, per, and $d_{\chi}^H$ on ${\cal C}_{n}$}\label{C}
\end{figure}

Let $A=(a_{ij})\in M_n$ and $k$ be a positive integer. Denote by $A^{[k]}$ the $k$-fold Hadamard product of $A$, that is, $A^{[k]}=(a_{ij}^k)$. If $A$ is a correlation matrix, then $\lim_{k\rightarrow \infty} A^{[k]}$ is a correlation matrix having entries 1, $-1$, and 0. It is easy to see that $\per \big (A^{[k]}\big ) \leq \per A$ if $k$ is large enough.
 If $A \in {\cal C}_{n}$ is such a matrix that has an maximizer $M_A$ for which $\per (A\circ M_A)>\per A$, then there must exist a
positive integer $r$ such that $\per (A\circ M_A^{[r]})\leq \per A$.
Drury's $16\times 16$ example shows that $\per (A\circ A)> \per A$ for some real correlation matrices $A$.

\begin{enumerate}
 \item[]{\bf Question 6:}
Characterize all the matrices $A\in {\cal C}_{n}$ that have the Hadamard-\indent compression property, that is,
$\per (A\circ X)\leq \per A $ for all $X\in {\cal C}_{n}$.
These are exactly the matrices having the property
$\per (A\circ M_A)\leq \per A$.

 \item[]{\bf Question 7:}
  Other than the determinant, are there any other generalized \indent matrix functions $d_{\chi}^H$ on ${\cal C}_{n}$ that are Hadamard-dilation?
  Namely, for all $n\times n$ correlation matrices $A$ and $B$, $d_{\chi}^H(A\circ B)\geq \max\{ d_{\chi}^H(A), d_{\chi}^H(B)\}.$

  Note that such a matrix function is closely related to the subgroup $H$ and character. If $H=\{e\}$ (the Hadamard matrix function),  then the inequality holds true.

 \item[]{\bf Question 8:} Characterize all maximizers of a given matrix $A\in {\cal C}_{n}$.
\end{enumerate}

Moreover, some other fragmentary results are seen.  It is known \cite[Result~5]{ZFZ89} that for   $n\times n$ correlation matrices $A$ and $B$,
$\per (A\circ B) \leq \per C_t$,
where $C_t$ is the $n\times n$  correlation matrix with all off-diagonal entries equal to $t$, where  $t=1-\l_{\min}(A),$ $ t=1-\l_{\min}(B),$ or $ t=1-\l_{\min}(A\circ B).$
It is also known  \cite[Result 4]{ZFZ89} that if $A$ is nonnegative entrywise, i.e., $a_{ij}\geq 0$, and $B$ is positive semidefinite, then
$|\per (A\circ B)|\leq \per (A) b_{11}\cdots b_{nn}.$

\section{Chollet conjecture is open}

In view of the classical result (of Oppenheim 1930)
 $\det (A\circ B)\geq \det A \, \det B$ for $A, B\geq 0$,
Chollet 1982 \cite{Cho82per}
 proposed an interesting problem concerning
the permanent of the Hadamard product of positive semidefinite
matrices:
\medskip

\noindent
{\bf Chollet conjecture 1982 \cite{Cho82per}.}
{\em  If $A\geq 0$ and $ B\geq 0$ are $n\times
n$ matrices,
then
\begin{equation}\label{Cho}
\per (A \circ B)\leq \per A\, \per B
\end{equation}
}

Chollet showed that the inequality holds if and only if it is true when $B=\bar{A}$:
\begin{equation}\label{Cho2}
\per (A\circ \bar{A})\leq (\per A)^2
\end{equation}

(\ref{Cho2}) is immediate from (\ref{Cho}) with $B=\bar{A}$. (\ref{Cho2}) implies  (\ref{Cho}) as shown in \cite{Cho82per}:
\begin{eqnarray}
0\leq \per (A\circ B) &  = & \sum_{\s \in S_n}\prod_{t=1}^n a_{t\s (t)} b_{t\s (t)} \nonumber \\
  & \leq   & \sum_{\s \in S_n}\prod_{t=1}^n |a_{t\s (t)}|\,| b_{t\s (t)}| \nonumber \\
  & \leq & \left ( \sum_{\s \in S_n}\prod_{t=1}^n |a_{t\s (t)}|^2\right )^{\frac12}\left ( \sum_{\s \in S_n}\prod_{t=1}^n |b_{t\s (t)}|^2\right )^{\frac12} \nonumber \\
  & = & \left ( \per (A\circ \bar{A})\right )^{\frac12} \left ( \per (B\circ \bar{B})\right )^{\frac12} \nonumber \\
  & \leq & \per A \, \per B
  \end{eqnarray}

With Drury's example in \cite{Dru16},
$\per (A\circ B)=\per (A\circ \bar{A})=\frac{6185}{128}<50<45^2=(\per A)^2$. Moreover, if $A=(a_{ij})$ is a correlation matrix, then all $|a_{ij}|\leq 1$. Thus,
$\per ( A\circ \bar{A})\leq \per (|A|)$, where $|A|=(|a_{ij}|)$. Note that
 $\per A\leq \per (|A|).$

Chollet conjecture is weaker than the Bapat \& Sunder conjecture because
$\per B\geq \prod_{i=1}^n b_{ii}$.
Special cases 
of (\ref{Cho}) 
  were settled in the 1980s.
  Gregorac and Hentzel 1987 \cite{GH87} showed
 by elementary methods that the inequality is true for the case of $2\times 2$ and $3\times 3$ matrices, which are also immediate consequences of the stronger results of  Bapat and Sunder 1986 \cite{BS86}. Grone and Merris 1987 \cite{GM87} and also
 Marcus and Sandy 1988 \cite{MarSan88}
 discussed  the Chollet conjecture and the Bapat \& Sunder conjecture.
  Related works are seen in Zhang 1989 \cite{ZFZ89} and Beasley
  2000 \cite{Bea00} (on some infinite classes of matrices).
 As the Bapat \& Sunder conjecture
  (\ref{BS85}) is settled in the negative recently,  Chollet conjecture stands out and is appealing as an interesting open question.
The real case of the Chollet conjecture is the following.

\begin{enumerate}
 \item[]{\bf Question 9:}
 Let $A=(a_{ij})$ be  a  positive semidefinite matrix of real entries. Does the following inequality hold?
$$\per \left ((a^2_{ij})\right ) 
 \leq (\per A)^2$$
\end{enumerate}

\section{A few more open conjectures on permanent}

Below are a few conjectures and research problems that virtually have no progress made in the past years. They are easily stated and understood, but appear to be forgotten. We bring them up again and hope to get renewed attention.
\medskip

\noindent
{\bf Marcus per-in-per conjecture 1965 \cite{MM65Per}} 
{\em  Let $A$ be a positive semidefinite $mn\times mn$ matrix partitioned as $A=(A_{ij})$
 in which each $A_{ij}$ is $n\times n$. Let $P$ be the $m\times m$ matrix whose
 $(i, j)$ entry is $\per A_{ij}$. Then
 $$\per A \geq \per P$$
 If $A$ is positive definite, then equality occurs if and only if all $A_{ij}=0$, $i\neq j$.
}
\medskip

 The case of $m=2$ is proved by Lieb 1966 \cite{Lieb66} and Djokovi\v{c} 1969 \cite{Djo69}. 
 A weaker inequality that $\per A\geq (\per P) /n!$ is shown by Pate 1981 \cite{Pate81}. 
 Pate 1982 \cite{Pate82} also shows that for each $m$, Marcus per-in-per conjecture is true for $n$ sufficiently large. 
 This may be compared with Thompson's det-in-det theorem 1961 \cite{Thm61} which asserts that
 the determinant of $A$ is dominated by the determinant of $D$ whose $(i,j)$ entry is $\det A_{ij}$, that is, $\det A\leq \det D$.

\medskip
\noindent
{\bf Marcus-Minc max-per-$U$ problem 1965 \cite{MM65Per}} 
{\em Let $A$ be an $n\times n$ positive semidefinite matrix. Find the maximum value of $\per (U^*AU)$ for all
$n\times n$ unitary matrices $U$.
}
\medskip

Permanent is a continuous function on the compact set of unitary matrices. The minimum and maximum values are achievable. Some upper bounds in terms of eigenvalues and trace of the matrix $A$ have been obtained. For instance,
$$\per (U^*AU)\leq \frac1n \sum_{i=1}^n\l^n_i$$
where $\l_i$ are the eigenvalues of the $n\times n$ positive semidefinite matrix $A$
(see, e.g., \cite{MarMin64}). 
It was conjectured that the maximum of $\per (U^*AU)$ is attained when all the main diagonal entries of $U^*AU$ are all equal (see \cite[p.\,132]{Minc87}). However, this is false; see a counterexample in \cite{Drew89b} or \cite[p.\,331]{CW05}. The problem is studied with partial solutions
by Drew and Johnson 1989 \cite{Drew89a} and Grone et al 1986 \cite{Grone86}.

\medskip
\noindent
{\bf Foregger per-$k$ conjecture 1978 \cite[p.\,157]{Minc78}.}
{\em
Let $n$ be a given positive integer. Then there exists a positive integer $k$ (solely depending on $n$) such that
for all $n\times n$ doubly stochastic matrices $A$
$$\per (A^k)\leq \per A$$
}

Chang 1984 \cite{Cha84} and 1990 \cite{Cha90} confirmed  the conjecture under certain conditions. For instance, if $\frac12 <\per A <1$, then $\per (A^m) < \per A$ for any integer $m\geq 2$. No other progress is seen.


\medskip
\noindent
{\bf Bapat-Sunder per-max conjecture 1986 \cite[Conjecture 3]{BS86}.}
{\em
Let $A$ be an $n\times n$ positive semidefinite matrix. Denote by $A(i,j)$ the $(n-1)\times (n-1)$ submatrix of $A$ obtained  by deleting the $i$th row and $j$th column of $A$. Then
$\per A$ is the maximum eigenvalue of the matrix  $\big (a_{ij}\per A(i,j)\big )$.
}
\medskip

This statement is weaker than the POT conjecture which is true for $n\leq 3$. Hence,
the Bapat-Sunder per-max conjecture holds true when the size of the matrix is no more than 3. Other than this, nothing else is known.

\medskip
For  nonnegative $A\in M_n$ and $B\in M_m$,  Brualdi 1966 \cite{Bru66} showed that
 $$(\per A)^m (\per B)^n \leq \per (A\otimes B)$$

 Does this also hold for positive semidefinite matrices?

\medskip
\noindent
{\bf Liang-So-Zhang per-tensor product conjecture 1992 \cite{ZLS94}.}
{\em  Let $A$ and $B$ be $n\times n$ and $m\times m$ positive semidefinite matrices, respectively. Then
 $$(\per A)^m (\per B)^n \leq \per (A\otimes B)$$
}

By Brualdi's result, the conjecture holds true for positive semidefinite matrices of nonnegative entries.
In particular, if $A$ and $B$ are square matrices of 1s,
it reduces to the known inequality $(n!)^m(m!)^n\leq (nm)!$.

Liang-So-Zhang  per-tensor product conjecture, if proven to be true, would imply a result of Ando 1981 \cite{Ando81}:
 $$\max\{ (n!)^{-m},\; (m!)^{-n}\}(\per A)^m (\per B)^n \leq \per (A\otimes B)$$
 which is stronger than a result of Marcus 1966 \cite{Mar66}:
 $$(n!)^{-m}(m!)^{-n}(\per A)^m (\per B)^n \leq \per (A\otimes B)$$

A special case of Liang-So-Zhang per-tensor product conjecture is the conjecture of Pate 1984 \cite{Pate84}:
If $M\in M_m$ is positive semidefinite and $J_k$ is the
all one matrix of size $k$, then
\begin{equation}\label{PateCon}
(k!)^m (\per M)^k\leq \per (J_k\otimes M)
\end{equation}
which has been confirmed for $k=2$ by Pate 1984 \cite{Pate84} and for $m=2$ \cite{ZLS94}.
$$2^m \left ( \per M\right )^2\leq \per \left (\begin{array}{cc}M & M \\
M & M
\end{array} \right ) $$

It would be interesting for one to consider the case where $M$ is $3\times 3$ in (\ref{PateCon}) and the case where $J_k$ is $3\times 3$ in (\ref{PateCon}).

\medskip
In light of the Chollet conjecture $\per (A\circ B)\leq \per A\per B$, we have

\medskip
\noindent
{\bf Liang-So-Zhang Hadamard-Kronecker product conjecture 1992 \cite{ZLS94}.}
{\em  Let $A$ and $B$ be $n\times n$ positive semidefinite matrices. Then
 $$\left (\per (A\circ B)\right )^n \leq \per (A\otimes B)$$}

The above conjecture is proved  \cite[Theorem 2]{ZLS94}
 for any positive semidefinite $A\in M_n$ and special $B=(b_{ij})\in M_n$ with all diagonal entries equal (to $a$, say) and
  all the entries above the main diagonal equal (to $b_{ij}=b$, say, for all $i<j$).

  \medskip
\noindent
{\bf Drury permanent conjecture 2016 \cite{Drury20160530}.}
{\em  Let $A$ be an $n\times n$ positive semidefinite matrix and $B$ be the submatrix of $A$ obtained by deleting the first row and column of $A$.  Let $B_{jk}$ be the submatrix of $B$
obtained by deleting the $j$th row and $k$th column of $B$. Then
\begin{equation}\label{Drury1}
(a_{11} \per B)^2 + \left ( \sum_{k=2}^n |a_{1k}|^2\per (B_{kk}) \right )^2  \leq (\per A)^2
\end{equation}
}
Two immediate consequences of (\ref{Drury1}) are $a_{11} \per B\leq \per A$ (well-known) and
\begin{equation}\label{Drury2}
\sum_{k=2}^n |a_{1k}|^2\per (B_{kk})\leq   \per A
\end{equation}
(\ref{Drury2}) itself is an interesting open problem. Chollet conjecture follows from (\ref{Drury1}).

\section*{Acknowledgements}
I am grateful to R. Bapat, R. Brualdi, G.-S. Cheon, S. Drury, and V. Shchesnovich
 for valuable comments. 
 I am especially indebted to Prof. Drury for his help with the counterexamples.
 The project was partially supported by National Natural Science Foundation of China (No.~11571220) through Shanghai University.

\end{document}